\documentclass{LMCS}

\usepackage{amsbsy}
\usepackage{enumerate}
\usepackage{bussproofs}
\usepackage{verbatim}
\usepackage{hyperref}

\newcommand{\ie}{i.e.}

\makeatletter
\def\verbatim{\small \@verbatim \frenchspacing\@vobeyspaces \@xverbatim}

\makeatother

\def\doframeit#1{\vbox{%
  \hrule height\fboxrule
    \hbox{%
      \vrule width\fboxrule \kern\fboxsep
      \vbox{\kern\fboxvsep #1\kern\fboxvsep }%
      \kern\fboxsep \vrule width\fboxrule }%
    \hrule height\fboxrule }}

\def\frameit{\smallskip \advance \linewidth by -7.5pt \setbox0=\vbox \bgroup
\strut \ignorespaces }

\def\endframeit{\ifhmode \par \nointerlineskip \fi \egroup
\doframeit{\box0}}

\newdimen \fboxvsep
\fboxvsep=\fboxsep
\advance \fboxsep by -3.5pt

\newcommand{\colspace}{.8cm}

\newcommand{\BNFDef}{::=}
\newcommand{\sep}{\; | \;}

\newcommand{\lam}[2]{\lambda {#1} . {#2}}
\newcommand{\ap}{\,}
\newcommand{\pOne}[1]{\pi_1 \, {#1}}
\newcommand{\pTwo}[1]{\pi_2 \, {#1}}
\newcommand{\pair}[2]{\langle {#1} , {#2} \rangle}

\newcommand{\pI}[1]{\pi_i \, {#1}}

\newcommand{\M}{M}
\newcommand{\N}{N}
\renewcommand{\P}{P}
\newcommand{\Q}{Q}

\newcommand{\FV}[1]{\mathrm{FV}(#1)}

\newcommand{\subst}[3]{{#1}[{#2} \mathrel{:=} {#3}]}

\newcommand{\rel}{\vartriangleright}

\newcommand{\genR}{\mathcal{R}}

\newcommand{\SP}{\protect{\text{\textsc{sp}}}}
\newcommand{\EP}{\text{\textsc{fp}}}

\newcommand{\SPl}{\text{\textsc{SP}}}
\newcommand{\EPl}{\text{\textsc{FP}}}

\newcommand{\etaExp}{\overline{\eta}}
\newcommand{\SPExp}{\overline{\SP}}

\newcommand{\main}{\EPl}

\newcommand{\BESP}{\beta \eta \SPl}
\newcommand{\BP}{\mathrm{R}}
\newcommand{\E}{\mathrm{E}}

\newcommand{\red}{\longrightarrow}
\newcommand{\redLeft}{\longleftarrow}

\newcommand{\redBP}{\red_{\BP}}
\newcommand{\redE}{\red_{\E}}
\newcommand{\redBESP}{\red_{\BESP}}
\newcommand{\redMain}{\red_{\main}}

\newcommand{\redELeft}{\redLeft_{\E}}

\newcommand{\SredBP}{\red_{\BP}^*}
\newcommand{\SredBE}{\red_{\beta \eta}^*}
\newcommand{\SredE}{\red_{\E}^*}

\newcommand{\SredMain}{\red_{\main}^*}

\newcommand{\parRed}{\Longrightarrow}
\newcommand{\parLeft}{\Longleftarrow}

\newcommand{\parE}{\parRed_{\E}}
\newcommand{\parBP}{\parRed_{\BP}}

\newcommand{\parELeft}{\parLeft_{\E}}

\newcommand{\equals}{=}

\newcommand{\eqBE}{\equals_{\beta \eta}}
\newcommand{\eqBESP}{\equals_{\beta \eta \SPl}}
\newcommand{\eqMain}{\equals_{\main}}

\newcommand{\theory}[1]{\lambda_{#1}}

\newcommand{\piSym}{ $\pi$-symmetric }

\newcommand{\erase}[1]{\lvert {#1} \rvert}

\theoremstyle{definition}
\newtheorem*{rem*}{Remark}
\newtheorem*{rems*}{Remarks}

\def\doi{2 (2:1) 2006}
\lmcsheading%
{\doi}
{14}
{}
{}
{Nov.~23, 2005}
{Mar.~10, 2006}
{}

\begin{document}

\title[Extensional Lambda Calculus with Surjective Pairing]
{Extending the Extensional Lambda Calculus\\ with Surjective
  Pairing\\ is Conservative}

\author[K.~St{\o}vring]{Kristian St{\o}vring}
\address{BRICS, Department of Computer Science, $\!$University of Aarhus\\
IT-parken, Aabogade~34, DK-8200 Aarhus N, Denmark}
\email{kss@brics.dk}


\begin{abstract}
  We answer Klop and de~Vrijer's question whether
  adding surjective-pairing axioms to the extensional lambda calculus
  yields a conservative extension.
  The answer is positive.
  As a byproduct we obtain a
  ``syntactic'' proof that the extensional lambda calculus with
  surjective pairing is consistent.
\end{abstract}

\maketitle
\vskip-\bigskipamount

\section{Introduction}
The theory $\theory{\BESP}$ is obtained from the untyped
extensional lambda calculus
$\theory{\beta\eta}$~\cite{Barendregt:84}
by adding three \emph{surjective-pairing} axioms:
\[
  \begin{array}[c]{crcll}
  (\pi_1) & \pOne{\pair{\M}{\N}} &=& \M &\\
  (\pi_2) & \pTwo{\pair{\M}{\N}} &=& \N &\\
  (\SP) & \pair{\pOne{\M}}{\pTwo{\M}} &=& \M &\\
  \end{array}
\]
These axioms are said to define a surjective pairing since the axiom
$(\SP)$ implies that every term is equal to a pair.

A $\lambda$-term is called
\emph{pure} if it does not contain any of the new constructs $\pi_i$
and $\pair{\cdot}{\cdot}$.  In this article we give a positive answer
to the following question, asked by Klop and de~Vrijer in
1989~\cite{Klop-de-Vrijer:IAC1989,de-Vrijer:LICS1989} and featured as
Problem 5 in the original RTA list of open
problems~\cite{Dershowitz-al:RTA1991}:
\begin{quote}
  Suppose that $\M$ and $\N$ are pure $\lambda$-terms.  Does $\/ \M
  \eqBESP \N$
  imply that $\M \eqBE \N$?
\end{quote}
In other words, we show that the theory $\theory{\BESP}$ is a
\emph{conservative extension} of the theory $\theory{\beta
  \eta}$.
As a byproduct we obtain a proof of consistency of $\theory{\BESP}$
that uses purely syntactic methods.\footnote{The author only knows of one
other such ``syntactic'' consistency proof for $\theory{\BESP}$,
namely one based on recent work on operationally defined
bisimulations~\cite{Lassen:LICS2006}.}
\vfill
\subsection{Background of the problem}
The two perhaps most obvious attempts at showing conservativity of
$\theory{\BESP}$ fail because of two
negative results:
no surjective-pairing function (that is, no pairing function satisfying the
three axioms above) is definable in the lambda~calculus~\cite{Barendregt:ZMLGM1974}, and
the standard reduction relation for the lambda calculus with surjective
pairing is not confluent~\cite{Klop:PhD-foo}.  Both results were shown for
the \emph{extensional} lambda calculus as well.

Klop~\cite{Klop:PhD-foo} and Klop and
de~Vrijer~\cite{Klop-de-Vrijer:IAC1989} have considered a number of
properties of the (non-extensional) lambda calculus with surjective
pairing, $\lambda_{\beta\SPl}$, which would have trivially followed
from confluence of the standard reduction relation. In particular,
de~Vrijer has shown that $\theory{\beta\SPl}$ is a conservative extension
of the lambda calculus~\cite{de-Vrijer:LICS1989}.  This result
motivated the question answered here: whether
surjective pairing also conservatively extends the \emph{extensional}
lambda calculus.

The proof of conservativity by de~Vrijer is furthermore the first
known ``syntactic'' consistency proof for $\theory{\beta \SPl}$. One of
Scott's model-theoretic consistency proofs for
$\theory{\beta\eta}$~\cite{Scott:CACM77-2} can be easily adapted to
show consistency of $\theory{\BESP}$ (and hence also $\theory{\beta \SPl}$)
as well.

The theory $\theory{\BESP}$ has also been investigated from a
categorical point of view.
If $\mathcal{C}$ is a cartesian closed
category with an object $U$ such that
\[U \cong U \times U \cong U \rightarrow U,\]
then there are various ways of interpreting $\lambda$-terms as morphisms
of $\mathcal{C}$~\cite{Barendregt:84,Lambek-Scott:86}.
Moreover, every
extension of the theory $\theory{\BESP}$ is the theory of a model
arising in this way~\cite{Lambek-Scott:86,Scott:Curry1980}.

\subsection{Formalization}
The author has formalized and verified the proof of the conservativity
result using the Twelf system~\cite{Pfenning-Schuermann:CADE1999}.
The formalized proof additionally serves as an implementation of a
procedure transforming a formal derivation of $\M \eqBESP N$ into a
formal derivation of $M \eqBE N$ (for pure terms $\M$ and $\N$).  It
is available from
\[ \text{\texttt{http://purl.oclc.org/net/kss/eta-SP}} \]
The formalized statement of the main result is presented in
Appendix~\ref{sec:formalized}.

\section{Background and notation}
The reader is assumed to be familiar with basic properties of the
untyped lambda calculus, as presented for example in the first three
chapters of Barendregt's book~\cite{Barendregt:84}.

The syntax of $\lambda$-terms is extended with constructs
for pairing and projection:
\[
     \M \; \BNFDef \; x \sep \lam{x}{\M} \sep \M \ap \M
      \sep \pair{\M}{\M} \sep \pOne{\M} \sep \pTwo{\M}
\]
(where $x$ ranges over an infinite set of variables).
The \emph{pure terms} are the usual $\lambda$-terms, \ie,
terms with no occurrences of $\pi_i$ or $\pair{\cdot}{\cdot}$.
The set of free variables of a term $\M$ is denoted $\FV{\M}$.
We follow practice and
identify $\alpha$-equivalent terms.

We use the following notation and definitions for relations on
$\lambda$-terms: For any binary relation $\rel_\genR$ on $\lambda$-terms,
$\red_\genR$ denotes the compatible closure of $\rel_\genR$ as defined in
Figure~\ref{fig:compatible}.  The relation $\red_\genR$ is called a
\emph{reduction relation}.  The reflexive--transitive closure of
$\red_\genR$ is written $\red_\genR^*$, and the reflexive-transitive-symmetric
closure of $\red_\genR$ is written $=_\genR$.
We write $\theory{\genR}$ for the
equational theory of $\lambda$-terms corresponding to $=_\genR$, i.e.,
$\theory{\genR}$ is the set of formal equations ``$\M = \N$'' such that
$\M =_\genR \N$.

\begin{figure}[htbp]
  \centering
  \begin{frameit}
\[
    \begin{array}[c]{c@{\quad\quad}c}
     % Base case
     \AxiomC{$\M \rel_\genR \M'$}
     \UnaryInfC{$\M \red_\genR \M'$}
     \DisplayProof &
    % Lambda
    \AxiomC{$\M \red_\genR \M'$}
    \UnaryInfC{$\lam{x}{\M} \red_\genR \lam{x}{\M'}$}
    \DisplayProof\\[\colspace]
    % App-1
    \AxiomC{$\M \red_\genR \M'$}
    \UnaryInfC{$\M \ap \N \red_\genR \M' \ap \N$}
    \DisplayProof &
    % App-2
    \AxiomC{$\N \red_\genR \N'$}
    \UnaryInfC{$\M \ap \N \red_\genR \M \ap \N'$}
    \DisplayProof\\[\colspace]
    % Pair-1
    \AxiomC{$\M \red_\genR \M'$}
    \UnaryInfC{$\pair{\M}{\N} \red_\genR \pair{\M'}{\N}$}
    \DisplayProof &
    % Pair-2
    \AxiomC{$\N \red_\genR \N'$}
    \UnaryInfC{$\pair{\M}{\N} \red_\genR \pair{\M}{\N'}$}
    \DisplayProof\\[\colspace]
    % Proj-1
    \AxiomC{$\M \red_\genR \M'$}
    \UnaryInfC{$\pOne{\M} \red_\genR \pOne{\M'}$}
    \DisplayProof &
    % Proj-2
    \AxiomC{$\M \red_\genR \M'$}
    \UnaryInfC{$\pTwo{\M} \red_\genR \pTwo{\M'}$}
    \DisplayProof
    \end{array}
\]
  \caption{The compatible closure of $\rel_\genR$.}
  \label{fig:compatible}
  \end{frameit}
\end{figure}
The relation $\rel_{\BESP}$ is defined by the
axioms in Figure~\ref{fig:BESP}.  This relation generates a
reduction relation $\redBESP$ and an equality relation $\eqBESP$.  The
extensional lambda calculus with surjective pairing is defined as the theory
$\theory{\BESP}$.

\begin{figure}[htbp]
  \centering
\begin{frameit}
\[
  \begin{array}[c]{crcll}
  (\beta) & (\lam{x}{\M}) \ap \N &\rel_{\BESP}& \subst{\M}{x}{\N} &\\
  (\eta) &  \lam{x}{\M} \ap x &\rel_{\BESP}& \M & \text{(if $x \notin \FV{\M}$)} \\
  (\pi_1) & \pOne{\pair{\M}{\N}} &\rel_{\BESP}& \M &\\
  (\pi_2) & \pTwo{\pair{\M}{\N}} &\rel_{\BESP}& \N &\\
  (\SP) & \pair{\pOne{\M}}{\pTwo{\M}} &\rel_{\BESP}& \M &\\
  \end{array}
\]
\caption{The relation $\rel_{\beta \eta \textsc{SP}}$.}
\label{fig:BESP}
\end{frameit}
\end{figure}

\section{Overview of the proof}
The relation $\red_{\BESP}$ 
is the standard reduction relation generating $=_{\BESP}$.
This reduction relation is, however, not
confluent~\cite{Curien-Hardin:JFP1994}~\cite[p.~216]{Klop:PhD-foo}; its
confluence would immediately imply the main result, namely that
$\theory{\BESP}$ is conservative over
$\theory{\beta\eta}$.\footnote{The non-confluent reduction relation
  considered by Klop~\cite{Klop:PhD-foo} is slightly different from
  $\red_{\BESP}$.  It is simple to construct
  a counter-example to confluence similar to Klop's.}

In this article we instead consider a
further extension $\theory{\main}$ of $\theory{\BESP}$ and show that
$\theory{\main}$ is conservative over $\theory{\beta\eta}$.
Since $\theory{\main}$ is an extension of $\theory{\BESP}$,
the main result follows.
The proof is structured in the following way:
\begin{itemize}
\item In Section~\ref{sec:confluence} we present the extension
  $\theory{\main}$ of $\theory{\BESP}$ and show that it is generated
  by a confluent reduction relation $\redMain$.  In the relation
  $\redMain$ the orientation of the axioms $(\eta)$ and $(\SP)$ is
  reversed; in other words, the extensionality axioms are oriented as
  \emph{expansion} axioms (see, e.g., the work by Jay and
  Ghani~\cite{Jay-Ghani:JFP95}).
\item In Section~\ref{sec:conservative} we show that $\theory{\main}$
  is conservative over $\theory{\beta \eta}$ on pure $\lambda$-terms.
  This result does not immediately follow from confluence of $\redMain$
  since $\redMain$ contains $(\SP)$ oriented as an expansion axiom.
\end{itemize}

\section{An extension of the theory $\theory{\beta \eta \textsc{SP}}$}
\label{sec:confluence}
We first present the extension $\theory{\main}$ of $\theory{\BESP}$;
the name $\EP$ is intended to be a mnemonic for ``functional
pairing''.
The relation $\rel_{\main}$ is defined by the axioms in
Figure~\ref{fig:BEP}.  This relation generates the theory
$\theory{\main}$ and the reduction relation~$\redMain$.
For convenience, we refer to the axioms
$(\delta \pi)$, $(\pi_1 \lambda)$, and $(\pi_2 \lambda)$ as the
\emph{commutation} axioms;  intuitively, these axioms express how a
function behaves as a pair and vice versa.
As discussed
above, the axioms $(\etaExp)$ and $(\SPExp)$ %% in $\redMain$
are oriented as expansion axioms.

The theory $\theory{\main}$ it not new, although is does not appear to
have been explicitly named before.  Axioms similar to the
commutation axioms $(\delta \pi)$, $(\pi_1 \lambda)$, and $(\pi_2
\lambda)$ were first considered in work on 
products and
lists in the lambda calculus~\cite{Revesz:TCS1993}
and in work on categorical
combinators~\cite{Revesz:FI1995}: adding the surjective-pairing axiom
$(\SP)$ to R\'ev\'esz's theory $\lambda_p$ gives the theory
$\theory{\main}$, except for a minor syntactic difference.  Durfee
gave a model for the full theory $\theory{\main}$~\cite{Durfee:MSc}
(see the remark below).  Axioms equivalent to the commutation axioms
play an indirect, but important, role in recent work on solvability for $\lambda$-terms
with pairs~\cite{Lassen:LICS2006}.

The reduction relation $\redMain$
(with its combination of commutation axioms and expansion axioms)
appears to be new.

\begin{figure}[htbp]
\begin{frameit}
\[
\begin{array}[c]{crcll}
  (\beta) & (\lam{x}{\M}) \ap \N &\rel_{\main}& \subst{\M}{x}{\N} &\\
  (\etaExp) & \M &\rel_{\main}& \lam{x}{\M} \ap x & \text{(if $x \notin \FV{\M}$)} \\
  (\pi_1) & \pOne{\pair{\M}{\N}} &\rel_{\main}& \M &\\
  (\pi_2) & \pTwo{\pair{\M}{\N}} &\rel_{\main}& \N &\\
  (\SPExp) & \M &\rel_{\main}& \pair{\pOne{\M}}{\pTwo{\M}} &\\[.2cm]
  (\delta\pi)     & \pair{\M}{\N} \ap \P &\rel_{\main}& \pair{\M \ap \P}{\N \ap \P} &\\
  (\pi_1 \lambda)      & \pOne{(\lam{x}{\M})} &\rel_{\main}& \lam{x}{\pOne{\M}} &\\
  (\pi_2 \lambda)      & \pTwo{(\lam{x}{\M})} &\rel_{\main}& \lam{x}{\pTwo{\M}} &\\
\end{array}
\]
\caption{The relation $\rel_{\main}$.}
\label{fig:BEP}
\end{frameit}\centering

\end{figure}

\begin{rem*}
  In this article, the theory $\theory{\main}$ and the associated
  reduction relation $\redMain$ are used to prove a specific result
  about a different theory. However, $\theory{\main}$ and $\redMain$
  can be justified semantically and syntactically:
  \begin{itemize}
  \item From the point of view of semantics: The original model of
    $\theory{\BESP}$~\cite{Lambek-Scott:86,Scott:CACM77-2}
    is
    also a model of $\theory{\main}$~\cite{Durfee:MSc}. Indeed, let
    $U$ and $V$ be complete partial orders such that $V \cong V \times
    V$ and $U \cong [U \rightarrow V]$.
Then by calculations valid in any cartesian
    closed category~\cite{Scott:CACM77-2}, $U \cong U \times U
    \cong [U \rightarrow U]$,
and one can verify 
that the standard interpretation\footnote{See
      also Exercise 18.4.19 in Barendregt's
      book~\cite{Barendregt:84}.} of $\lambda$-terms as elements of
    $U$ gives rise to a model of $\theory{\main}$.
    
    As an aside, if $U$ is an arbitrary complete partial order satisfying
    that $U \cong U \times U \cong [U \rightarrow U]$, then the
    standard interpretation using these isomorphisms makes $U$ a model
    of (at least) $\theory{\BESP}$.  Taking $U = V$ in the above
    construction now gives an alternative pair of isomorphisms, and
    hence an alternative interpretation of $\lambda$-terms, resulting
    in a model of $\theory{\main}$.
    
  \item From the point of view of term rewriting: In the
    \emph{simply-typed} lambda calculus, term constructs can be
    proof-theoretically classified as either \emph{introduction forms}
    ($\lam{x}{\M}$ and $\pair{\M}{\N}$) or \emph{elimination forms}
    ($\M \ap \N$ and $\pI{\M}$), using the Curry-Howard
    isomorphism~\cite{Barendregt:92-foo}.  The simply-typed counterparts
    of the axioms $(\beta)$, $(\pi_1)$, and $(\pi_2)$ of
    Figure~\ref{fig:BEP} then imply that, when constructing a term
    bottom-up, ``an introduction form followed by an elimination form
    is a redex.''  This property is preserved in the untyped reduction
    relation $\redMain$ by virtue of the commutation axioms $(\delta \pi)$,
    $(\pi_1 \lambda)$ and $(\pi_2 \lambda)$.
  \end{itemize}
\end{rem*}

In the rest of this section we prove that $\redMain$ is confluent.
For that purpose we describe~$\redMain$ as the union of two relations:
a part $\redE$ generated from the $\eta$/$\SP$-expansion axioms
$(\etaExp)$ and $(\SPExp)$, and an ``extensionality-free'' part
$\redBP$ generated from all the remaining axioms.
\begin{itemize}
\item In Section~\ref{sec:confluence-BP} we
  show that the extensionality-free part $\redBP$ is confluent.
\item In Section~\ref{sec:eta-postponement} we review the
  well-known fact that $\eta$/$\SP$-expansion $\redE$ is confluent,
  and then show that $\redE$ commutes with $\redBP$:
  if $\N_1 \redELeft^* \M \SredBP
  \N_2$, then there is a $\P$ such that $\N_1 \SredBP \P \redELeft^*
  \N_2$.
\end{itemize}
We conclude by
the Hindley--Rosen Lemma~\cite[p.~64]{Barendregt:84} that the union
$\mathord{\redMain} = \mathord{\redBP} \cup \mathord{\redE}$ is
confluent.  Earlier, van~Oostrom used a similar approach to prove confluence of
$\eta$-expansion (together with $\beta$-reduction)
in the pure lambda calculus~\cite{van-Oostrom:TCS1997}.

From a technical point of view, the proof that $\redE$
commutes with $\redBP$ is the novel part of the confluence proof: the
commutation proof highlights the role of the axioms $(\delta\pi)$,
$(\pi_1 \lambda)$, and $(\pi_2 \lambda)$.

\subsection{Confluence of an extensionality-free subrelation}
\label{sec:confluence-BP}
The relation $\rel_{\BP}$ is defined by all the axioms of
$\rel_{\main}$ except $(\etaExp)$ and $(\SPExp)$; for convenience the
remaining axioms are shown in Figure~\ref{fig:BP}.  The relation
$\rel_{\BP}$ generates the reduction relation $\redBP$.
\begin{figure}[htbp]
\begin{frameit}
\[
\begin{array}[c]{crcll}
  (\beta) & (\lam{x}{\M}) \ap \N &\rel_{\BP}& \subst{\M}{x}{\N} &\\
  (\pi_1) & \pOne{\pair{\M}{\N}} &\rel_{\BP}& \M &\\
  (\pi_2) & \pTwo{\pair{\M}{\N}} &\rel_{\BP}& \N &\\
  (\delta\pi)       & \pair{\M}{\N} \ap \P &\rel_{\BP}& \pair{\M \ap \P}{\N \ap \P} &\\
  (\pi_1 \lambda)   & \pOne{(\lam{x}{\M})} &\rel_{\BP}& \lam{x}{\pOne{\M}} &\\
  (\pi_2 \lambda)   & \pTwo{(\lam{x}{\M})} &\rel_{\BP}& \lam{x}{\pTwo{\M}} &\\
\end{array}
\]
\caption{The relation $\rel_{\BP}$.}
\label{fig:BP}
\end{frameit}\centering

\end{figure}

We now aim to prove that $\redBP$ is confluent.
In fact, this follows from general higher-order rewriting theory,
since $\redBP$ can be formulated as an orthogonal pattern higher-order
rewriting system~\cite{Nipkow:TLCA1993,vanRaamsdonk:PhD},
and such systems are
confluent~\cite{Nipkow:TLCA1993}.
However, in order to keep the presentation self-contained,
we give a direct confluence proof.
This direct proof, which follows the method of the
Tait/\mbox{Martin-L{\"o}f} proof of confluence of $\beta$-reduction
~\cite[p.~60]{Barendregt:84}, can be viewed as
a specialized version of Nipkow's confluence
proof~\cite{Nipkow:TLCA1993}.
\vfill\eject

First, define a parallel~\cite{Takahashi:IaC1995} reduction relation
$\parBP$, shown in Figure~\ref{fig:parBP}.\footnote{The notion that
$\parBP$ is the parallel reduction relation generated from the axioms
of $\rel_{\BP}$ can be made precise~\cite[Section 4]{Nipkow:TLCA1993}.}

\begin{figure}[htbp]
  \centering
  \begin{frameit}
\[
    \begin{array}[c]{c@{\quad\quad}c}
      \multicolumn{2}{c}{%
       \AxiomC{$\M \parBP \M'$}
       \AxiomC{$\N \parBP \N'$}
       \BinaryInfC{$(\lam{x}{\M}) \ap \N \parBP \subst{\M'}{x}{\N'}$}
       \DisplayProof}\\[\colspace]
      %% P1
      \AxiomC{$\M \parBP \M'$}
      \UnaryInfC{$\pOne{\pair{\M}{\N}} \parBP \M'$}
      \DisplayProof &
      %% P2
      \AxiomC{$\N \parBP \N'$}
      \UnaryInfC{$\pTwo{\pair{\M}{\N}} \parBP \N'$}
      \DisplayProof \\[\colspace]
      %% App-pair
      \multicolumn{2}{c}{
      \AxiomC{$\M \parBP \M'$}
      \AxiomC{$\N \parBP \N'$}
      \AxiomC{$\P \parBP \P'$}
      \TrinaryInfC{$\pair{\M}{\N} \ap \P \parBP
        \pair{\M' \ap \P'}{\N' \ap \P'}$}
      \DisplayProof} \\[\colspace]
      %% proj1-lam
      \AxiomC{$\M \parBP \M'$}
      \UnaryInfC{$\pOne{(\lam{x}{\M})} \parBP \lam{x}{\pOne{\M'}}$}
    \DisplayProof &
      %% proj2-lam
      \AxiomC{$\M \parBP \M'$}
      \UnaryInfC{$\pTwo{(\lam{x}{\M})} \parBP \lam{x}{\pTwo{\M'}}$}
    \DisplayProof \\[\colspace]
      %% Refl
      \AxiomC{}
      \UnaryInfC{$\M \parBP \M$}
      \DisplayProof &
      %% Lam
      \AxiomC{$\M \parBP \M'$}
      \UnaryInfC{$\lam{x}{\M}  \parBP \lam{x}{\M'}$}
    \DisplayProof\\[\colspace]
      %% App
      \AxiomC{$\M \parBP \M'$}
      \AxiomC{$\N \parBP \N'$}
      \BinaryInfC{$\M \ap \N  \parBP \M' \ap \N'$}
    \DisplayProof &
      %% Pair
      \AxiomC{$\M \parBP \M'$}
      \AxiomC{$\N \parBP \N'$}
      \BinaryInfC{$\pair{\M}{\N}  \parBP \pair{\M'}{\N'}$}
      \DisplayProof
      \\[\colspace]
      %% Proj1
      \AxiomC{$\M \parBP \M'$}
      \UnaryInfC{$\pOne{\M}  \parBP \pOne{\M'}$}
    \DisplayProof &
      %% Proj2
      \AxiomC{$\M \parBP \M'$}
      \UnaryInfC{$\pTwo{\M}  \parBP \pTwo{\M'}$}
    \DisplayProof\\[\colspace]
    \end{array}
\]
  \caption{Parallel $\BP$-reduction $\parBP$.}
  \label{fig:parBP}
  \end{frameit}
\end{figure}
\begin{prop}
  \label{prop:parBP-properties}\mbox{}
  \begin{enumerate}[(i)]
  \item
$\SredBP \; = \; \parBP^*$.
  \item If\/ $\M \parBP \M'$ and $\N \parBP \N'$, then
  $\subst{\M}{x}{\N} \parBP \subst{\M'}{x}{\N'}$.
\item If\/ $\M \SredBP \M'$ and $\N \SredBP \N'$, then $\subst{\M}{x}{\N}
  \SredBP \subst{\M'}{x}{\N'}$.
  \end{enumerate}
\end{prop}
\proof
  Standard~\cite[p.~60]{Barendregt:84}.  Part (iii) follows from the
  first two parts and will be used in the next section.
\qed
\begin{prop}
  The relation $\parBP$ satisfies the diamond property: if\/
  $M \parBP N_1$ and $M \parBP N_2$, then there is a $P$ such that
  $N_1 \parBP P$ and $N_2 \parBP P$.
\end{prop}
\proof
  By induction on the derivations of $M \parBP N_1$ and $M \parBP N_2$
  according to the rules in Figure~\ref{fig:parBP}.  Many of the cases
  are well-known from the proof of confluence of \mbox{$\beta$-reduction}.
  There are no interesting new cases (which is another way of saying
  that $\redBP$ can naturally be defined as an \emph{orthogonal}
  higher-order term rewriting system).
\qed

\begin{cor}
\label{prop:confluent-BP}
  The relation $\redBP$ is confluent.
\end{cor}

\subsection{The relation $\redBP$ commutes with 
  $\eta$/\raise 2 pt\hbox{$_{\textrm{SP}}$}-expansion}
%\subsection{The relation $\redBP$ commutes with $\eta$/$\SP$-expansion}
\label{sec:eta-postponement}
We define the relation $\rel_{\E}$ by the axioms $(\etaExp)$ and
$(\SPExp)$, for convenience shown in Figure~\ref{fig:E}.  This
relation generates the $\eta$/$\SP$-expansion relation $\redE$.
\begin{figure}[htbp]
\begin{frameit}
\[
\begin{array}[c]{crcll}
  (\etaExp) & \M &\rel_{\E}& \lam{x}{\M} \ap x & \text{(if $x \notin \FV{\M}$)} \\
  (\SPExp) &  \M &\rel_{\E}& \pair{\pOne{\M}}{\pTwo{\M}} &\\[.2cm]
\end{array}
\]
\caption{The relation $\rel_{\E}$.}
\label{fig:E}
\end{frameit}\centering

\end{figure}

The purpose of this section is to show that $\redE$ commutes with
$\redBP$, that is, if $\N_1 \redELeft^* \M \SredBP \N_2$, then there is a $\P$
such that $\N_1 \SredBP \P \redELeft^* \N_2$.  Before proceeding with
the proof of commutation, we consider some of the critical
pairs~\cite{Nipkow:TLCA1993}
between $\redE$ and~$\redBP$.  The first two cases are well-known:
\begin{enumerate}
\item $(\lam{x}{(\lam{x}{\M}) \ap x}) \ap \N \redLeft_{\etaExp} (\lam{x}{\M}) \ap \N
  \red_{\beta} \subst{\M}{x}{\N}$.\\[.05cm]
  Solution:
$
(\lam{x}{(\lam{x}{\M}) \ap x}) \ap \N
  \red_{\beta} (\lam{x}{\M}) \ap \N
  \red_{\beta} \subst{\M}{x}{\N}.
$\vskip.1cm
\item $\pI{\pair{\pOne{\pair{\M_1}{\M_2}}}{\pTwo{\pair{\M_1}{\M_2}}}}
  \redLeft_{\SPExp} \, \; \pI{\pair{\M_1}{\M_2}} \red_{\pi_i}
  \M_i$.\\[.05cm]
  Solution:
  $
  \pI{\pair{\pOne{\pair{\M_1}{\M_2}}}{\pTwo{\pair{\M_1}{\M_2}}}}
  \red_{\pi_i} \pI{\pair{\M_1}{\M_2}} \red_{\pi_i}
  \M_i.
  $
\end{enumerate}
On the other hand, to resolve the next
two kinds of critical pairs, one
needs the commutation axioms $(\delta \pi)$, $(\pi_1 \lambda)$, and $(\pi_2
\lambda)$:
\begin{enumerate}
\setcounter{enumi}{2}
\item $\pOne{(\lam{x}{\pair{\M_1}{\M_2} \ap x})} \redLeft_{\etaExp} \;
  \pI{\pair{\M_1}{\M_2}} \red_{\pi_i} \M_i$.\\[.05cm]
  Solution:
\[
  \begin{array}[c]{rll}
\pI{(\lam{x}{\pair{\M_1}{\M_2} \ap x})} & \red_{\delta\pi} &
      \pI{(\lam{x}{\pair{\M_1 \ap x}{\M_2 \ap x}})}\\
      & \red_{\pi_i \lambda} & \lam{x}{\pI{\pair{\M_1 \ap x}{\M_2 \ap
      x}}}\\
      & \red_{\pi_i} & \lam{x}{\M_i \ap x}\\
      & \redLeft_{\etaExp} & \M_i.
  \end{array}
\]
\item $\pair{\pOne{(\lam{x}{\M})}}{\pTwo{(\lam{x}{\M})}} \ap \N
  \redLeft_{\SPExp} (\lam{x}{\M}) \ap \N \red_{\beta}
  \subst{\M}{x}{\N}$.\\[.05cm]
Solution:
\[
  \begin{array}[c]{rll}
    \pair{\pOne{(\lam{x}{\M})}}{\pTwo{(\lam{x}{\M})}} \ap \N
    & \red_{\pi_1 \lambda , \pi_2 \lambda}^* &
    \pair{\lam{x}{\pOne{\M}}}{\lam{x}{\pTwo{\M}}} \ap \N\\
    & \red_{\delta\pi} & \pair{(\lam{x}{\pOne{\M}}) \ap \N}
                              {(\lam{x}{\pTwo{\M}}) \ap \N}\\
    & \red_{\beta}^* & \pair{\pOne{(\subst{\M}{x}{\N})}}
                            {\pTwo{(\subst{\M}{x}{\N})}}\\
    & \redLeft_{\SPExp} & \subst{\M}{x}{\N}.\\
  \end{array}
\]
\end{enumerate}
These are all the kinds of critical pairs between $\redE$ and $\redBP$
in which the $\BP$-step uses one of the axioms $(\beta)$, $(\pi_1)$,
or $(\pi_2)$.  The cases where the $\BP$-step is one of the remaining
axioms can be resolved similarly to the simple cases $1$ \mbox{and $2$}.

We now turn to the actual proof of commutation.  Define
a parallel $\eta$/$\SP$-expansion relation
$\parE$~\cite{Jay-Ghani:JFP95,Takahashi:IaC1995} by the rules in
Figure~\ref{fig:parE}.
\begin{figure}[htbp]
\begin{frameit}
  \centering
\[
    \begin{array}[c]{c@{\quad\quad}c}
    %% Eta
      \AxiomC{$\M \parE \M'$}
      \RightLabel{$\;$ ($x \notin \FV{\M}$)}
      \UnaryInfC{$\M  \parE  \lam{x}{\M' \ap x}$}
      \DisplayProof &
    %% SP
      \AxiomC{$\M \parE \M'$}
      \UnaryInfC{$\M  \parE  \pair{\pOne{\M'}}{\pTwo{\M'}}$}
      \DisplayProof \\[1.2cm]
      %% Refl
      \AxiomC{}
      \UnaryInfC{$\M \parE \M$}
      \DisplayProof &
      %% Lam
      \AxiomC{$\M \parE \M'$}
      \UnaryInfC{$\lam{x}{\M}  \parE \lam{x}{\M'}$}
    \DisplayProof\\[\colspace]
      %% App
      \AxiomC{$\M \parE \M'$}
      \AxiomC{$\N \parE \N'$}
      \BinaryInfC{$\M \ap \N  \parE \M' \ap \N'$}
    \DisplayProof &
      %% Pair
      \AxiomC{$\M \parE \M'$}
      \AxiomC{$\N \parE \N'$}
      \BinaryInfC{$\pair{\M}{\N}  \parE \pair{\M'}{\N'}$}
      \DisplayProof
      \\[\colspace]
      %% Proj1
      \AxiomC{$\M \parE \M'$}
      \UnaryInfC{$\pOne{\M}  \parE \pOne{\M'}$}
    \DisplayProof &
      %% Proj2
      \AxiomC{$\M \parE \M'$}
      \UnaryInfC{$\pTwo{\M}  \parE \pTwo{\M'}$}
    \DisplayProof\\[\colspace]
    \end{array}
\]
  \caption{Parallel $\eta$/$\SP$-expansion $\parE$.}
  \label{fig:parE}
\end{frameit}
\end{figure}

First, some simple facts about parallel $\eta$/$\SP$-expansion:
\begin{prop}
  \label{prop:parE-properties}
  \label{prop:eta-ud-confluent}\hfill
  \begin{enumerate}[(i)]
  \item $\redE^* \; = \; \parE^*$.
  \item $\redE$ is confluent.
  \item If $\M \parE \M'$ and $\N \parE \N'$, then
  $\subst{\M}{x}{\N} \parE \subst{\M'}{x}{\N'}$.
  \end{enumerate}
\end{prop}
\proof
  Standard~\cite{Jay-Ghani:JFP95}.
The confluence of $\redE$ follows from
 the diamond property of $\parE$.
\qed

We now aim to prove that if $\N_1 \parELeft \M \redBP \N_2$, then
there exists a $\P$ such that $\N_1 \SredBP \P \parELeft \N_2$.  Consider
for example the case
\[
\N \ap \Q \parELeft (\lam{x}{\M}) \ap \Q \red_{\beta} \subst{\M}{x}{\Q}
\]
where $\N \parELeft \lam{x}{M}$. Then $\N$ results from $\lam{x}{\M}$
by a number of $\eta$/$\SP$-expansions, and in order to show
commutation we intuitively need to iterate cases $1$ and $4$ of the
critical pair calculations shown in the beginning of this section.  Similar
examples exist for the other axioms of $\redBP$.  The properties which
are needed are shown in the following two lemmas:
\begin{lem}
  \label{lem:parE-first}
  If\/ $\lam{x}{\M} \parE \N$, then
  \begin{enumerate}[(i)]
  \item there is a $\P$ such that $\N \ap x \SredBP P \parELeft \M$, and
  \item there is a $\Q$ such that for $i \in \{1,2\}$: $\pI{\N} \SredBP
    \lam{x}{\pI{\Q}}$ and $\M \parE \Q$.
  \end{enumerate}
\end{lem}
\proof
By induction on the definition of $\lam{x}{\M} \parE \N$.
\qed

\begin{lem}
\label{lem:parE-last}
 If\/ $\pair{\M_1}{\M_2} \parE \N$, then
  \begin{enumerate}[(i)]
  \item for $i \in \{1,2\}$ there is a $\P_i$ such that $\pI{\N} \SredBP
    P_i \parELeft \M_i$, and
  \item there are $\Q_1$, $\Q_2$ such that $\N \ap x \SredBP
    \pair{\Q_1 \ap x}{\Q_2 \ap x}$ and $\M_1 \parE \Q_1$ and
    $\M_2 \parE \Q_2$.
\end{enumerate}
\end{lem}
\proof
By induction on the definition of $\pair{\M_1}{\M_2} \parE \N$.
\qed

We now prove the main lemma needed in the commutation proof:
\begin{lem}
\label{lem:eta-postponement}
  If\/ $\N \parELeft \M \redBP \M'$, then there is a $\P$ such that
  $\N \SredBP \P \parELeft \M'$.
\end{lem}
\proof
  Induction on the definition of $\M \parE \N$, using
  Lemmas~\ref{lem:parE-first} and \ref{lem:parE-last}.  We
  show some illustrative cases.
\begin{enumerate}[C{a}se 1:]
\item $\pair{\pOne{\N'}}{\pTwo{\N'}} \parELeft \M \redBP \M'$ where $\N'
  \parELeft \M$.  By the induction hypothesis there is a $\P'$ such that 
  $\N' \SredBP \P' \parELeft \M'$.  Then
\[
\pair{\pOne{\N'}}{\pTwo{\N'}} \SredBP \pair{\pOne{\P'}}{\pTwo{\P'}}
\parELeft \M'
\]
so choose $\P = \pair{\pOne{\P'}}{\pTwo{\P'}}$.
\item $\N_1 \ap \N_2 \parELeft (\lam{x}{\M_1}) \ap \M_2 \redBP
  \subst{\M_1}{x}{\M_2}$ where $\N_1 \parELeft \lam{x}{\M_1}$ and where
  \mbox{$\N_2 \parELeft \M_2$}.
By Lemma~\ref{lem:parE-first}(i) there is a
  $\P'$ such that $\N_1 \ap x \SredBP \P' \parELeft \M_1$.
  It is easy to see from the definition of $\mathrel{\parE}$ that $x$
  is not free in~$\N_1$.
Therefore, 
by
  Propositions~\ref{prop:parBP-properties}
  and~\ref{prop:parE-properties}, $\N_1 \ap \N_2 \SredBP
  \subst{\P'}{x}{\N_2} \parELeft \subst{\M_1}{x}{\M_2}$, so choose $\P
  = \subst{\P'}{x}{\N_2}$. \qed
\end{enumerate}

\begin{lem}
\label{lem:eta-postponement-star}\mbox{}
\begin{enumerate}[(i)]
\item   If\/ $\N \parELeft \M \SredBP \M'$, then there is a $\P$ such that
  $\N \SredBP \P \parELeft \M'$.
\item If\/ $\N \parELeft^* \M \SredBP \M'$, then there is a $\P$ such
  that $\N \SredBP \P \parELeft^* \M'$.
\end{enumerate}
\end{lem}
\proof\mbox{}
  \begin{enumerate}[(i)]
  \item   By induction on the length of the reduction sequence $\M \SredBP
  \M'$, using Lemma~\ref{lem:eta-postponement}.
  \item By induction on the length of the reduction sequence $\M \parE^*
    \N$, using Part (i). \qed
\end{enumerate}

Now, by Proposition~\ref{prop:parE-properties}(i), $\SredE \!\; = \; \parE^*$.
Therefore Lemma~\ref{lem:eta-postponement-star}(ii) implies
that the relations $\redE$ and $\redBP$ commute:
\begin{prop}
  \label{prop:eta-post-1}
  If\/ $\N \redELeft^* \M \SredBP \M'$, then there is a $\P$ such that
  $\N \SredBP \P \redELeft^* \M'$.
\end{prop}

\subsection{Confluence of $\redMain$}
\label{sec:confluence-main}
We now use the results of
Sections~\ref{sec:confluence-BP}~and~\ref{sec:eta-postponement} to
prove the main result of Section~\ref{sec:confluence}:
\begin{prop}
\label{prop:confluent-main}
  The relation $\redMain$ is confluent.
\end{prop}
\proof
  Proposition~\ref{prop:confluent-BP} states that $\redBP$ is
  confluent, Proposition~\ref{prop:eta-ud-confluent}(ii) states that
  $\redE$ is confluent, and Proposition~\ref{prop:eta-post-1} states
  that $\redBP$ commutes with $\redE$.  By the Hindley--Rosen
  Lemma~\cite[p.~64]{Barendregt:84}, the relation $\mathord{\redMain}  =
  \mathord{\redBP} \cup \mathord{\redE}$ is confluent. 
\qed

\begin{cor}[Church--Rosser property]
  \label{cor:CR}
  If\/ $\M \eqMain \N$, then there is a $\P$ such that $\M \SredMain \P$
  and $\N \SredMain \P$.
\end{cor}
\proof
  Follows from confluence of $\redMain$~\cite[p.~54]{Barendregt:84}.
\qed

\begin{rems*}\hfill
\begin{enumerate}[(i)]
\item Orienting the axioms $(\SPExp)$ and $(\etaExp)$ of $\redMain$ as
  \emph{contraction} axioms does not give rise to a confluent
  reduction relation: with these axioms we would have the reductions $\lam{x}{x} \redLeft_{\main}
  \pair{\pOne{(\lam{x}{x})}}{\pTwo{(\lam{x}{x})}} \SredMain
  \pair{\lam{x}{\pOne{x}}}{\lam{x}{\pTwo{x}}}$, but the two terms
  $\lam{x}{x}$ and $\pair{\lam{x}{\pOne{x}}}{\lam{x}{\pTwo{x}}}$ would
  be normal forms.
\item The commutation axioms of $\theory{\main}$ depend on the fact
  that the calculus is untyped, such that, intuitively, every function
  is also a pair and vice versa.  A different line of work concerns
  reduction relations in \emph{typed} calculi, with product and unit types,
  containing $(\SP)$ oriented as a contraction
  axiom~\cite{Curien-DiCosmo:JFP1996}.
\end{enumerate}
\end{rems*}

\section{Main result}
\label{sec:conservative}
We are now almost in a position to prove the main result: Suppose $\M$
and $\N$ are pure $\lambda$-terms such that $\M \eqBESP \N$.  Then $\M
\eqMain \N$, and by the Church--Rosser \mbox{property}
\mbox{(Corollary~\ref{cor:CR})} there is a $\P$ such that $\M \SredMain \P$
and $\N \SredMain \P$.  However, since~$\redMain$ contains
$\SP$-\emph{expansion}, we cannot immediately conclude that $P$ is a
pure \mbox{$\lambda$-term} with $\M
\SredBE \P$ and $\N \SredBE \P$.

\begin{defi}%%[$\pi$-erasure]
  The \emph{$\pi$-erasure} of a $\lambda$-term $\M$ is the pure
  $\lambda$-term $\erase{\M}$ defined inductively as follows:
\[
\begin{array}[c]{rcl}
  \erase{x} &=& x\\
  \erase{\M \ap \N} &=& \erase{\M} \ap \erase{\N}\\
  \erase{\lam{x}{\M}} &=& \lam{x}{\erase{\M}}\\
  \erase{\pair{\M}{\N}} &=& \erase{\M}\\
  \erase{\pOne{\M}} &=& \erase{\M}\\
  \erase{\pTwo{\M}} &=& \erase{\M}\\
\end{array}
\]
\end{defi}
We could just as well have defined $\erase{\pair{\M}{\N}}$ as
$\erase{\N}$, since we are only interested in $\erase{\P}$ when $\P$ is
\emph{$\pi$-symmetric}:
\begin{defi}
  A $\lambda$-term $\M$ is \emph{$\pi$-symmetric} if for every subterm
  of $\M$ of the form $\pair{\P}{\Q}$,
the $\pi$-erasures of $P$ and $Q$ are $\beta\eta$-equivalent:
$\erase{\P} \eqBE \erase{\Q}$.
\end{defi}
In particular, every pure $\lambda$-term is $\pi$-symmetric.

\begin{prop}
  \label{prop:pisym-subst}\hfill
  \begin{enumerate}[(i)]
  \item $\erase{\subst{\M}{x}{\N}} = \subst{\erase{\M}}{x}{\erase{\N}
  \hskip 1pt}$
  \item If\/ $\M$ and $\N$ are $\pi$-symmetric, then $\subst{\M}{x}{\N}$
  is $\pi$-symmetric.
  \end{enumerate}
\end{prop}
\proof
  By induction on $\M$.
\qed

\begin{prop}
  \label{prop:pisym-reduction}
  If\/ $M$ is \piSym and $M \redMain N$, then
  \begin{enumerate}[(i)]
  \item $\erase{\M} \eqBE \erase{\N}$, and
  \item $N$ is $\pi$-symmetric.
  \end{enumerate}
\end{prop}
\proof
  By induction on the definition of $M \redMain N$, using
  Proposition~\ref{prop:pisym-subst}.
\qed

Now we are ready to prove that $\theory{\main}$ is a conservative extension of
$\theory{\beta\eta}$:
\begin{thm}
\label{thm:main-conservative}
  Let $\M, \N$ be pure $\lambda$-terms.  If\/ $\M \eqMain N$, then $M \eqBE N$.
\end{thm}
\proof
  Suppose $\M$ and $\N$ are pure $\lambda$-terms such that $\M \eqMain
  \N$. By the Church--Rosser property (Corollary~\ref{cor:CR}) there is
  a $\P$ such that $\M \SredMain \P$ and $\N \SredMain \P$.  Since
  $\M$ and $\N$ are pure, they are in particular $\pi$-symmetric; it
  follows from Proposition~\ref{prop:pisym-reduction} that $\P$
  is $\pi$-symmetric and that $\erase{\M} \eqBE \erase{\P} \eqBE
  \erase{N}$.  Hence
$\M = \erase{\M}
\eqBE \erase{N} = \N$.  \qed

\begin{cor}
  The theory $\theory{\main}$ is consistent.
\end{cor}
\proof
  By Theorem~\ref{thm:main-conservative} and consistency of
  $\theory{\beta\eta}$~\cite[p.~67]{Barendregt:84}.
\qed
Finally we turn to the main result of this article:
\begin{thm}
  Let $\M, \N$ be pure $\lambda$-terms.  If\/ $\, M \eqBESP N$, then $M \eqBE N$.
\end{thm}
\proof
  By Theorem~\ref{thm:main-conservative} and the fact that
  $\theory{\main}$ is an extension of $\theory{\BESP}$.
\qed
We have also obtained a new---syntactic---proof of consistency of
$\theory{\BESP}$:
\begin{cor}
  The theory $\theory{\BESP}$ is consistent.
\end{cor}

\begin{rem*}
  \label{rem:beta-eta-D}
  The question of conservativity was originally formulated in a
  slightly different setting~\cite{Klop-de-Vrijer:IAC1989}: let $D$,
  $D_1$ and $D_2$ be three new constants, and add the following axioms
  to the pure $\lambda_{\beta\eta}$-calculus:
\[
\begin{array}[c]{rcl}
  D_1 \ap (D \ap \M \ap \N) &=_{\beta\eta D}& \M\\
  D_2 \ap (D \ap \M \ap \N) &=_{\beta\eta D}& \N\\
  D \ap (D_1 \ap \M) \ap (D_2 \ap \M) &=_{\beta\eta D} & \M
\end{array}
\]
To see that the resulting theory $\theory{\beta\eta D}$ is
conservative over $\theory{\beta\eta}$, one can simulate
$\theory{\beta \eta D}$ in $\theory{\BESP}$ by defining $D$
as $\lam{x}{\lam{y}{\pair{x}{y}}}$, $D_1$ as $\lam{x}{\pOne{x}}$, and
$D_2$ as $\lam{x}{\pTwo{x}}$.
\end{rem*}

\section{Related problems}

The conservativity proof presented here can be adapted to the
non-extensional case settled by de~Vrijer~\cite{de-Vrijer:LICS1989},
i.e., a minor modification gives an alternative proof that
$\theory{\beta\SPl}$ is conservative over the lambda calculus
$\theory{\beta}$.
To this end, one should simply remove the axiom $(\eta)$ from every
definition and proof.
The electronic, formalized version of the proof allows
for a straightforward verification that the modification is correct.

Another related problem posed by Klop and de~Vrijer is still open:
whether the reduction relation $\red_{\BESP}$
has the \emph{unique normal-form
  property}~\cite{Klop-de-Vrijer:IAC1989}.
The theory $\theory{\main}$ does not seem useful in solving that problem.

Meyer asked whether \emph{any} lambda theory can be conservatively
extended with surjective pairing~\cite{Dershowitz-al:RTA1991}.  That
problem also remains open.

\section*{Acknowledgements}
The author is grateful to Olivier Danvy, Andrzej Filinski,
and the anonymous referees for their insightful comments.
Vincent van~Oostrom pointed out a substantial simplification of the
confluence proof in \mbox{Section~\ref{sec:confluence}}.  Thanks are
also due to Pierre-Louis Curien and Soren Lassen for discussions on
this work, and to
Karl Crary for his lectures on LF and the Twelf system in the fall of
2004 at CMU.

The work described in this article is supported by BRICS (Basic
Research in Computer Science ({\tt {http://www.brics.dk}}), funded by
the Danish National Research Foundation).

\bibliographystyle{plain}

\begin{thebibliography}{10}

\bibitem{Barendregt:ZMLGM1974}
Henk Barendregt.
\newblock Pairing without conventional restraints.
\newblock {\em Z. Math. Logik Grundlag. Math.}, 20:289--306, 1974.

\bibitem{Barendregt:84}
Henk Barendregt.
\newblock {\em The Lambda Calculus: Its Syntax and Semantics}, volume 103 of
  {\em Studies in Logic and the Foundation of Mathematics}.
\newblock North-Holland, revised edition, 1984.

\bibitem{Barendregt:92-foo}
Henk Barendregt.
\newblock Lambda calculi with types.
\newblock In Samson Abramsky, Dov~M. Gabbay, and Thomas S.~E. Maibaum, editors,
  {\em Handbook of Logic in Computer Science, Vol. 2}, chapter~2, pages
  118--309. Oxford University Press, Oxford, 1992.

\bibitem{Curien-DiCosmo:JFP1996}
Pierre-Louis Curien and Roberto Di~Cosmo.
\newblock A confluent reduction system for the lambda-calculus with surjective
  pairing and terminal object.
\newblock {\em Journal of Functional Programming}, 6(2):299--327, 1996.

\bibitem{Curien-Hardin:JFP1994}
Pierre-Louis Curien and Th{\'e}r{\`e}se Hardin.
\newblock Yet yet a counterexample for $\lambda$+{SP}.
\newblock {\em Journal of Functional Programming}, 4(1):113--115, 1994.

\bibitem{Dershowitz-al:RTA1991}
Nachum Dershowitz, Jean-Pierre Jouannaud, and Jan~Willem Klop.
\newblock Open problems in rewriting.
\newblock In Ronald~V. Book, editor, {\em Rewriting Techniques and
  Applications, 4th International Conference, RTA-91}, volume 488 of {\em
  Lecture Notes in Computer Science}, pages 445--456. Springer-Verlag, 1991.
\newblock The RTA list of open problems is currently maintained at
  \texttt{http://www.lsv.ens-cachan.fr/rtaloop/}.

\bibitem{Durfee:MSc}
Glenn Durfee.
\newblock A model for a list-oriented extension of the lambda calculus.
\newblock Master's thesis, School of Computer Science, Carnegie Mellon
  University, 1997.

\bibitem{Jay-Ghani:JFP95}
C.~Barry Jay and Neil Ghani.
\newblock The virtues of eta-expansion.
\newblock {\em Journal of Functional Programming}, 5(2):135--154, 1995.

\bibitem{Klop:PhD-foo}
Jan~Willem Klop.
\newblock {\em Combinatory Reduction Systems}.
\newblock Mathematical Centre Tracts 127. Mathematisch Centrum, Amsterdam,
  1980.

\bibitem{Klop-de-Vrijer:IAC1989}
Jan~Willem Klop and Roel de~Vrijer.
\newblock Unique normal forms for lambda calculus with surjective pairing.
\newblock {\em Information and Computation}, 80(2):97--113, 1989.

\bibitem{Lambek-Scott:86}
Joachim Lambek and Philip~J. Scott.
\newblock {\em Introduction to Higher Order Categorical Logic}, volume~7 of
  {\em Cambridge studies in advanced mathematics}.
\newblock Cambridge University Press, 1986.

\bibitem{Lassen:LICS2006}
Soren~B. Lassen.
\newblock Head normal form bisimulation for pairs and the
  \mbox{$\lambda\mu$-calculus}.
\newblock Manuscript, 2006.

\bibitem{Nipkow:TLCA1993}
Tobias Nipkow.
\newblock Orthogonal higher-order rewrite systems are confluent.
\newblock In Marc Bezem and Jan~Friso Groote, editors, {\em Typed Lambda
  Calculi and Applications, TLCA '93}, volume 664 of {\em Lecture Notes in
  Computer Science}, pages 306--317. Springer-Verlag, 1993.

\bibitem{van-Oostrom:TCS1997}
Vincent {\SortNoop{Oostrom}}van~Oostrom.
\newblock Developing developments.
\newblock {\em Theoretical Computer Science}, 175(1):159--181, 1997.

\bibitem{Pfenning:CR}
Frank Pfenning.
\newblock A proof of the {C}hurch--{R}osser theorem and its representation in a
  logical framework.
\newblock Technical Report CMU-CS-92-186, School of Computer Science, Carnegie
  Mellon University, 1992.

\bibitem{Pfenning-Schuermann:CADE1999}
Frank Pfenning and Carsten Sch{\"u}rmann.
\newblock System description: Twelf - a meta-logical framework for deductive
  systems.
\newblock In Harald Ganzinger, editor, {\em Automated Deduction---CADE-16, 16th
  International Conference on Automated Deduction}, volume 1632 of {\em Lecture
  Notes in Computer Science}, pages 202--206. Springer-Verlag, 1999.

\bibitem{vanRaamsdonk:PhD}
Femke {\SortNoop{Raamsdonk}}van~Raamsdonk.
\newblock {\em Confluence and Normalization for higher-order rewriting}.
\newblock PhD thesis, Vrije Universiteit Amsterdam, 1996.

\bibitem{Revesz:TCS1993}
Gy{\"o}rgy~E. R{\'e}v{\'e}sz.
\newblock A list-oriented extension of the lambda-calculus satisfying the
  {C}hurch-{R}osser theorem.
\newblock {\em Theoretical Computer Science}, 93:75--89, 1992.

\bibitem{Revesz:FI1995}
Gy{\"o}rgy~E. R{\'e}v{\'e}sz.
\newblock Categorical combinators with explicit products.
\newblock {\em Fundamenta Informaticae}, 22:153--166, 1995.

\bibitem{Scott:CACM77-2}
Dana~S. Scott.
\newblock Logic and programming languages.
\newblock {\em Communications of the ACM}, 20:634--641, 1977.

\bibitem{Scott:Curry1980}
Dana~S. Scott.
\newblock Relating theories of the lambda calculus.
\newblock In J.~P. Seldin and J.~R. Hindley, editors, {\em To H.B. Curry:
  Essays on Combinatory Logic, Lambda-Calculus and Formalism}, pages 403--450.
  Academic Press, 1980.

\bibitem{Takahashi:IaC1995}
Masako Takahashi.
\newblock Parallel reductions in $\lambda$-calculus.
\newblock {\em Information and Computation}, 118:120--127, 1995.

\bibitem{de-Vrijer:LICS1989}
Roel {\SortNoop{Vrijer}}de~Vrijer.
\newblock Extending the lambda calculus with surjective pairing is
  conservative.
\newblock In {\em Proceedings of the Fourth Annual {IEEE} Symposium on Logic in
  Computer Science}, pages 204--215, Pacific Grove, California, June 1989. IEEE
  Computer Society Press.

\end{thebibliography}
\newcommand{\SortNoop}[1]{}

\appendix
\section{Formalized statement of the main result}
\label{sec:formalized}
Below is the formalized statement of the conservativity theorem.  The
full formal proof consists of 2670 lines of Twelf code.  It was developed
using version 1.5R1 of the Twelf system.\footnote{The Twelf system can
be obtained from \texttt{http://www.cs.cmu.edu/$\sim$twelf/}} The
encoding technique is based on a formal proof of the
Church--Rosser theorem for $\beta$-reduction that is distributed
along with earlier versions of the Twelf system~\cite{Pfenning:CR}.

\begin{verbatim}
%%% Terms of the untyped lambda calculus with surjective pairing.
\end{verbatim}
\begin{verbatim}
term : type.
\end{verbatim}
\begin{verbatim}
@ : term -> term -> term.  %infix left 10 @.
lam : (term -> term) -> term.
p1 : term -> term.
p2 : term -> term.
pair : term -> term -> term.
\end{verbatim}
\begin{verbatim}
%freeze term.
\end{verbatim}
\vskip.1cm

\begin{verbatim}
%%% Lambda calculus with the extensionality rules eta and SP.
\end{verbatim}
\begin{verbatim}
==SP : term -> term -> type.  %infix none 5 ==SP.
\end{verbatim}
\begin{verbatim}
sp_beta : (lam F) @ N ==SP F N.
\end{verbatim}
\begin{verbatim}
sp_eta : lam ([x] M @ x) ==SP M.
\end{verbatim}
\begin{verbatim}
sp_proj1 : p1 (pair M N) ==SP M.
\end{verbatim}
\begin{verbatim}
sp_proj2 : p2 (pair M N) ==SP N.
\end{verbatim}
\begin{verbatim}
sp_SP : pair (p1 M) (p2 M) ==SP M.
\end{verbatim}
\begin{verbatim}
% Congruence rules.
\end{verbatim}
\begin{verbatim}
sp_refl : M ==SP M.
\end{verbatim}
\begin{verbatim}
sp_sym : M ==SP N -> N ==SP M.
\end{verbatim}
\begin{verbatim}
sp_trans : M ==SP N -> N ==SP P -> M ==SP P.
\end{verbatim}
\begin{verbatim}
sp_c-app : M @ N ==SP M' @ N'
            <- M ==SP M'
            <- N ==SP N'.
\end{verbatim}
\begin{verbatim}
sp_c-lam : lam F ==SP lam F'
            <- ({x} F x ==SP F' x).
\end{verbatim}
\begin{verbatim}
sp_c-p1 : p1 M ==SP p1 M'
           <- M ==SP M'.
\end{verbatim}
\begin{verbatim}
sp_c-p2 : p2 M ==SP p2 M'
           <- M ==SP M'.
\end{verbatim}
\begin{verbatim}
sp_c-pair : pair M N ==SP pair M' N'
             <- M ==SP M'
             <- N ==SP N'.
\end{verbatim}
\begin{verbatim}
%freeze ==SP.
\end{verbatim}
\vskip.1cm

\begin{verbatim}
%%% Pure lambda-terms, i.e., no "pair", "p1", or "p2".
\end{verbatim}
\begin{verbatim}
pterm : type.
\end{verbatim}
\begin{verbatim}
^ : pterm -> pterm -> pterm.  %infix left 10 ^.
lambda : (pterm -> pterm) -> pterm.
\end{verbatim}
\begin{verbatim}
%freeze pterm.
\end{verbatim}
\vskip.1cm

\begin{verbatim}
%%% Beta-eta equality on pure terms.
\end{verbatim}
\begin{verbatim}
==be : pterm -> pterm -> type.  %infix none 5 ==be.
\end{verbatim}
\begin{verbatim}
be_beta : (lambda F) ^ N ==be F N.
\end{verbatim}
\begin{verbatim}
be_eta : lambda ([x] M ^ x) ==be M.
\end{verbatim}
\begin{verbatim}
% Congruence rules.
\end{verbatim}
\begin{verbatim}
be_refl : M ==be M.
\end{verbatim}
\begin{verbatim}
be_sym : M ==be N -> N ==be M.
\end{verbatim}
\begin{verbatim}
be_trans : M ==be N -> N ==be P -> M ==be P.
\end{verbatim}
\begin{verbatim}
be_c-app : M ^ N ==be M' ^ N'
           <- M ==be M'
           <- N ==be N'.
\end{verbatim}
\begin{verbatim}
be_c-lam : lambda F ==be lambda F'
           <- ({x} F x ==be F' x).
\end{verbatim}
\begin{verbatim}
%freeze ==be.
\end{verbatim}
\vskip.1cm

\begin{verbatim}
%%% Injecting pure terms into the general terms.
\end{verbatim}
\begin{verbatim}
inject : pterm -> term -> type.
%mode inject +P -T.
\end{verbatim}
\begin{verbatim}
inj_app : inject (P1 ^ P2) (M1 @ M2)
           <- inject P1 M1
           <- inject P2 M2.
\end{verbatim}
\begin{verbatim}
inj_lam : inject (lambda P) (lam M)
           <- ({x} {y} inject x y -> inject (P x) (M y)).
\end{verbatim}
\begin{verbatim}
%freeze inject.
\end{verbatim}
\begin{verbatim}
%block inj : block {x : pterm} {y : term} {thm : inject x y}.
\end{verbatim}
\begin{verbatim}
%worlds (inj) (inject _ _).
%total P (inject P _).
\end{verbatim}
\vskip.1cm

\begin{verbatim}
%%% The main theorem: ==SP is conservative over ==be.
\end{verbatim}
\begin{verbatim}
conservative : inject M M' -> inject N N'
                           -> M' ==SP N'
                           -> M ==be N
                           -> type.
%mode conservative +I1 +I2 +E1 -E2.
\end{verbatim}
\begin{verbatim}
% [The proof is omitted.]
\end{verbatim}
\begin{verbatim}
%worlds () (conservative _ _ _ _).
%total I1 (conservative I1 _ _ _).
\end{verbatim}
\begin{verbatim}
% With empty "worlds", the main theorem is actually only shown
% for closed terms. (The generalization to open terms is more
% complicated to express, but it follows easily by
% lambda-abstracting every free variable.)
\end{verbatim}
\vskip-40 pt\phantom{.}

\end{document}